\documentclass[a4paper,12pt]{article}
\usepackage{amsmath,amssymb}

\newtheorem{thrm}{Theorem}
\newtheorem{cor}{Corollary}
\newtheorem{lem}{Lemma}
\newtheorem{prop}{Proposition}
\newtheorem{defn}{Definition}
\newtheorem{exm}{Example}

\newcommand{\NatS}{{{\cal{N}{\!\it at}}S}}
\newcommand{\up}[2]{\stackrel{#2}{#1}}
\newcommand{\toop}{\up{\to}{\mbox{\large .}}}

\begin{document}

\title{\bf 0-Cohomology of semigroups}

\author{\bf B.\,Novikov}
\date{\normalsize Department of Mathematics, University of Kharkov, Ukraine\\
e-mail: {\it boris.v.novikov@univer.kharkov.ua}}

\maketitle

\tableofcontents

\bigskip

\bigskip

In 1974 my teacher L.\,M.\,Gluskin has offered me to study
projective representations of semigroups. As one could expect, cohomology appeared in
this problem, however it ``slightly'' differed
from the Eilenberg--MacLane cohomology. I have named it
``0-cohomology'', have studied its properties insofar as this was
necessary for the initial problem, and thought that I will probably never meet
this notion again.

However, for the last 30 years I returned to 0-cohomology again and
again since I met problems, in which it appeared.

This article is a survey of 0-cohomology. The main attention
is devoted to applications. Necessary definitions and results from
Homological Algebra and Theory of Semigroups can be found in
\cite{c-e}, \cite{c-p}, and \cite{h-s}.

I have to note that there are a lot of various cohomology theories
adapted to solution of specific problems in semigroups
\cite{cla,gri,han,lau,lee,log,pac,str,wel}; see also my review
\cite{nov13}. So an invention of one more cohomology is not any
novelty. However, 0-cohomology seems attractive (at least for me!)
because it links semigroups with different branches of the algebra.

\section{Eilenberg--MacLane cohomology}\label{sec-1}

The definition of semigroup cohomology does not differ from group
cohomology \cite{c-e}: for a semigroup $S$ and a (left) $S$-module
$A$ we call by an {\it $n$-dimensional cohomology group} the group
$H^n(S,A)={\rm Ext}^n_{\mathbb{Z}S}(\mathbb{Z}, A)$ where
$\mathbb{Z}$ is considered as a trivial $\mathbb{Z}S$-module. We
will name this cohomology by {\it Eilenberg--MacLane cohomology} or
briefly {\it EM-cohomology}.

In what follows we will use another well-known definition. By
$C^n(S,A)$ the group of all $n$-place mappings $f:\underbrace{S
\times \ldots \times S}_{n\ \mbox{\rm \footnotesize times}} \to A$
(the group of $n$-dimensional cochains) is denoted; a coboundary
operator $\partial ^n:C^n(S,A)\to C^{n+1}(S,A)$ is defined as
follows:
\begin{eqnarray}\label{eq-1-1}
&&\partial^nf(x_1,\ldots ,x_{n+1})=x_1f(x_2,\ldots ,x_{n+1}) \nonumber\\
&&+\sum_{i=1}^n{(-1)^if(x_1,\ldots ,x_ix_{i+1},\ldots ,x_{n+1})} +
(-1)^{n+1}f(x_1,\ldots ,x_n)
\end{eqnarray}

 Then $\partial^{n}\partial^{n-1}=0$, i.\,e.
$$
{\rm Im}\partial^{n-1}=B^n(S,A)\ \mbox{(the group of}\ n
\mbox{-dimensional coboundaries)}
$$
$$
\subseteq {\rm Ker}\partial^n=Z^n(S,A)\ \mbox{(the group of}\ n
\mbox{-dimensional cocycles)}
$$
and cohomology groups are defined as $H^n(S,A)=Z^n(S,A)/B^n(S,A)$.

However, for cohomology of semigroups one does not manage to obtain
results which are comparable to theory of cohomology of groups. In
this section we give several results about cohomology of semigroups,
illustrating its specifics.

Here is a typical example: since a projective module over a
semigroup is not obliged to be projective over its subsemigroup, the
lemma of Shapiro \cite{bro}, which expresses cohomology of a
subgroup through cohomology of groups, does not hold for semigroups.
So a result, received in \cite{c-e}, is of interest:

\begin{thrm}\label{thrm-1-1}
Let $T$ be a subsemigroup of a group $G$. If $G$ is a group of
fractions for $T$ (i.\,e. each element from $G$ can be written in
the form $x^{-1}y$ for some $x,y\in T$), then homomorphisms
$i^n:H^n(G,A)\to H^n(T,A)$, induced by the embedding $i:T\to G$, are
bijective for any $G$-module $A$.
\end{thrm}

Let $I$ be an ideal of a semigroup $S$. What can one say about
homomorphisms $\varepsilon^n:H^n(S,A)\to H^n(I,A)$, induced by the
embedding $\varepsilon :I\to S$? It is easy to show that
$\varepsilon^0$ is always an isomorphism and $\varepsilon^1$ is a
monomorphism. Using technique of adjoint functors, Adams and Rieffel
\cite{a-r} proved

\begin{thrm}\label{thrm-1-2}
Let $I$ be a left ideal of a semigroup $S$, having an identity $e$.
Then for any $S$-module $A$ and for any $n\ge 0$
$$
H^n(S,A)\cong H^n(I,A)\cong H^n(I,eA).
$$
\end{thrm}

In particular, if $S$ contains zero then $H^n(S,A)=0$ for $n>0$.

In \cite{a-r} by help of Theorem \ref{thrm-1-2}  a sufficient
condition was obtained for an associative algebra over $\mathbb{R}$
be a semigroup algebra.

The connection between $H^n(S,A)$ and $H^n(I,A)$ becomes more close
when we take for $I$ the so-called Sushkevich kernel (the least
two-sided ideal). This situation was in detail considered by
W.\,Nico \cite{nic}. I will not formulate this result, but only note
that it implies the description of cohomology of completely simple
semigroups:

\begin{thrm}\label{thrm-1-3}
Let $S$ be a completely simple semigroup, $G$ its basic group, $e$ the
identity of $G$, $A$ an $S$-module. Then $H^n(S,A)\cong H^n(G,eA)$
for any $n\ge 3$.
\end{thrm}

A number of papers is devoted to study of cohomological dimension
of semigroups. {\it Cohomological dimension} (cd$S$) of a semigroup
$S$ is defined as the least integer $n$ such that $H^{n+1}(S,A)=0$
for any $S$-module $A$.

It is well-known \cite{c-e} that cohomological dimension of both a free
group and a free semigroup (monoid) does not exceed 1. In the case
of groups the converse statement is also true (this is the famous
Stallings--Swan theorem  \cite{bro}). For semigroups the situation
is evidently different: joining an extra zero to any semigroup makes trivial all
its cohomology groups, except 0-dimensional. So it is
naturally to confine ourselves to the class of semigroups with
cancellation. Mitchell \cite{mit} has shown that the free product of
a free group and a free monoid (which he called a partially free
monoid) has cohomological dimension 1. In the same paper he has
formulated a suggestion: if $S$ is a cancellative monoid with ${\rm
cd}S=1$ then $S$ is partially free.

In \cite{nov5} (and later in \cite{nov7} more generally) a
counter-example to the Mitchell suggestion was built and a
``weakened Mitchell conjecture'' was proposed: if $S$ is a
cancellative monoid with ${\rm cd}S=1$ then $S$ is embedded into a
free group. This suggestion was proved in \cite{nov12}. Probably, it
is difficult to get more exact information: it is shown in
\cite{nov8} that a semigroup antiisomorphic to the counter-example
from \cite{nov5} (and therefore also embeddable into a free group)
has cohomological dimension 2. A good answer is known only in the
commutative case: cohomological dimensions of all subsemigroups of
$\mathbb{Z}$ are 1 \cite{nov9}.

\section{Properties of 0-cohomology}\label{sec-2}

Before constructing 0-cohomology we have to define a suitable
Abelian category. In what follows $S$ is a semigroup with a zero.

\begin{defn}\label{defn-2-1}
A $0$-module over $S$ is an Abelian (additive) group $A$ equipped with a
multiplication \mbox{$(S\setminus 0)\times A\to A$}
satisfying the following
conditions for all $s,t\in S\setminus 0,\ a,b\in A$:
$$
s(a+b)=sa+sb,
$$
$$
st\ne 0 \Rightarrow s(ta)=(st)a.
$$
A morphism of $0$-modules is a homomorphism of Abelian groups
$\varphi:A\to B$ such that $\varphi(sa)=s\varphi(a)$ for $s\in
S\setminus 0,\ a\in A$.
\end{defn}

We will denote the obtained category of $0$-modules by ${\rm
Mod}_0\,S$. It is easy to see that for the semigroup $T^0=T\cup0$
with an extra zero the category ${\rm Mod}_0\,T^0$ is isomorphic to
the category ${\rm Mod}\,T$ of usual modules over $T$. It turns out
that in general case ${\rm Mod}_0\,S$  is also isomorphic to a
category of (usual) modules over some semigroup.

Denote by $\overline{S}$ the set of all finite sequences
$(x_1,\ldots,x_m)$ such that $x_i\in S\setminus 0$\, $(1\le i\le
m)$ and $x_ix_{i+1}=0$ \,$(1\le i<m)$; thus, all one-element
sequences, except $(0)$, contain in $\overline{S}$. Define a binary relation $\rho$ on
$\overline{S}$ via
$(x_1,\ldots,x_m)\rho (y_1,\ldots,y_n)$ if and only if one of the
following conditions is fulfilled:

\medskip

1) $m=n$ and there exists $i$ $(1\le i \le m-1)$ that

$x_i=y_iu$, $y_{i+1}=ux_{i+1}$ for some $u\in S$,

and $x_j=y_j$ for $j\ne i$, $j\ne i+1$;

\medskip

2) $m=n+1$ and there exists $i$ $(2\le i \le m-1)$ that

$x_i=uv$, $y_{i-1}=x_{i-1}u$, $y_i=vx_{i+1}$ for some $u,v\in S$,

$x_j=y_j$ for $1\le j\le i-2$,

and $x_j=y_{j-1}$ for $i+2\le j \le m$.

\medskip

Let $\overline{\rho}$ be the least equivalence containing $\rho$,
$\tilde{S}$ the quotient set $\overline{S}/\overline{\rho}$. The
image of $(x_1,\ldots,x_m) \in \overline{S}$ in $\tilde{S}$ will be
denoted by $[x_1,\ldots,x_m]$.

Define a multiplication on $\tilde{S}$:
$$
[x_1,\ldots,x_m][y_1,\ldots,y_n]=\left\{
\begin{array}{ll}
[x_1,\ldots,x_my_1,\ldots,y_n],& \mbox{\rm if \ }\ x_my_1\ne 0,\\

[x_1,\ldots,x_m,y_1,\ldots,y_n],& \mbox{\rm if \ }\ x_my_1=0.
\end {array} \right.
$$
Then $\tilde{S}$ becomes a semigroup, which is called {\it a gown}
of $S$.

Each 0-module over $S$ can be transformed into an (usual) module
over $\tilde{S}$ by
$$
[x_1,\ldots,x_m]a=x_1(\ldots(x_ma)\ldots)
$$
for $x_1,\ldots,x_m \in S\setminus 0$, $a\in A$. Hence we obtain

\begin{prop}\label{prop-2-1}
${\rm Mod}_0\,S \cong {\rm Mod}\,\tilde{S}$.
\end{prop}

\begin{cor}\label{cor-2-1}
The category ${\rm Mod}_0\,S$ is Abelian.
\end{cor}

Here are some simplest properties of the gown \cite{nov6}:

\medskip

1. If $S=T^0$ is a semigroup with an extra zero then $\tilde{S}\cong
S\setminus 0 =T$.

2. It follows from the definition of the relation $\rho$ that the
map $S\setminus 0\to\tilde{S}$, $x\to [x]$, is bijective.

3. The subset $J=\{[x_1,\ldots,x_m]\in \tilde{S} \mid m>1\}$ is an
ideal in $\tilde{S}$ and $\tilde{S}/J\cong S$.

\medskip

It is easy to find the gown if the semigroup $S$ is given by
defining relations. We will write $S=\langle a_1,\ldots,a_m \mid
P_i=Q_i, 1\le i \le n\rangle$ if $S$ is generated by elements
$a_1,\ldots,a_m$ and is defined by equalities
$P_1=Q_1,\ldots,P_n=Q_n$. If the value of a word $P_i$ (or the same,
of $Q_i$) in semigroup $S$ is 0 then the equality $P_i=Q_i$ is
called {\it zero}.

\begin{prop}\label{prop-2-2}
Let $S=\langle a_1,\ldots,a_m \mid P_i=Q_i, 1\le i \le n\rangle$ be
a semigroup with a zero, in which none of generating elements is 0.
If one deletes all zero defining relations then the obtained semigroup
will be isomorphic to the gown $\tilde{S}$.
\end{prop}

\begin{exm}\label{exm-2-1}
If $S$ is a semigroup with zero multiplication ($S^2=0$) then
$\tilde{S}$ is a free semigroup.
\end{exm}

Let now $A$ be a 0-module over $S$.

\begin{defn}\label{defn-2-2}
A partial $n$-place mapping from $S$ into $A$, defined on all
$n$-tuples $(x_1,\ldots,x_n)$ such that $x_1\cdots x_n\ne 0$, is called
a $n$-dimensional cochain. The group of $n$-dimensional cochains is
denoted by $C^n_0(S,A)$. A coboundary operator $\partial
^n:C^n_0(S,A)\to C^{n+1}_0(S,A)$ is given as above, by the formula
(\ref{eq-1-1}). Then $\partial^{n}\partial^{n-1}=0$ and
$0$-cohomology groups are defined as
$H^n_0(S,A)=Z^n_0(S,A)/B^n_0(S,A)$, where $Z^n_0(S,A)= {\rm
Ker}\partial^n$ is the group of $n$-dimensional 0-cocycles,
$B^n_0(S,A)={\rm Im}\partial^{n-1}$ is the group of $n$-dimensional
$0$-coboundaries.
\end{defn}

\begin{exm}\label{exm-2-2}
Let $S=T^0=T\cup0$ be a semigroup with an extra zero. Then
$\tilde{S}\cong T$ and one can easily check that $H^n_0(S,A) \cong
H^n(T,A)$\footnote {Here we consider $A$ both as a 0-module over $S$
and a module over $\tilde{S}$, which does not lead to misunderstanding
in this context.} while $H^n(S,A)=0$.
\end{exm}

This example shows that 0-cohomology is a generalization of
EM-cohomo\-logy. In view of Proposition \ref{prop-2-1}, in general
case it is naturally to compare the groups $H^n_0(S,A)$ and
$H^n(\tilde{S},A)$. We describe this comparison in more details.

Since the sequence of functors $\{H^n_0(S,\_)\}_{n\ge 0}$ is
connected by terminology of \cite{h-s} and
$\{H^n(\tilde{S},\_)\}_{n\ge 0}$ is a sequence of derived functors
in the isomorphic categories ${\rm Mod}_0\,S$ and ${\rm
Mod}\,\tilde{S}$ respectively, then an isomorphism
$$
\varepsilon^0:H^0(\tilde{S},A)={\rm Hom}_{{\rm
Mod}\,\tilde{S}}(\mathbb{Z}, A) \cong {\rm Hom}_{{\rm
Mod}_0\,S}(\mathbb{Z}, A)=H^0_0(S,A)
$$
induces group homomorphisms
$$
\varepsilon^n:H^n(\tilde{S},A) \to H^n_0(S,A)
$$
such that $\{\varepsilon^n\}_{n\ge 0}$ are morphisms of cohomology
functors.

The homomorphism $\varepsilon^n$ is described as follows. If $f\in
C^n(\tilde{S},A)$ then set $(\eta^nf)(x_1,\ldots,x_n)=
f([x_1],\ldots,[x_n])$ for $x_1\cdot\ldots\cdot x_n \ne 0$. Then
$\eta^n$ is a homomorphism from $C^n(\tilde{S},A)$ into
$C_0^n(S,A)$ and it induces the homomorphism $\varepsilon^n$.

Direct calculations prove
\begin{thrm}\label{thrm-2-1}
$\varepsilon^1$ is an isomorphism for any semigroup $S$.
\end{thrm}

Using the corresponding long exact sequence, from Theorem
\ref{thrm-2-1} we obtain
\begin{cor}\label{cor-2-2}
$\varepsilon^2$ is a monomorphism for any semigroup $S$.
\end{cor}

Generally speaking, the groups $H^2_0(S,A)$ and $H^2(\tilde{S},A)$
are not isomorphic:

\begin{exm}\label{exm-2-3}
Let a commutative semigroup $S$ consists of elements $u,v,w,0$ with
the multiplication
$$
u^2=v^2=uv=w, \ uw=vw=0.
$$
One can show that $H^2(\tilde{S},A)=0$ for any module $A$ over $S$.
On the other hand, if $A$ (considered now as a $0$-module) is not
one-element and $a\in A\setminus 0$, then $0$-cocycle $f$ defined
by the condition
$$
f(x,y)=\left\{
\begin{array}{l}
a \ \mbox{\rm for} \ x=y=u,\\ 0\ \mbox{\rm otherwise},
\end {array} \right.
$$
is not a $0$-coboundary and thus $H^2_0(S,A)\ne 0$.
\end{exm}

Apropos, this example shows that 0-cohomology is not a derived
functor unlike EM-cohomology. Indeed, according to Proposition
\ref{prop-2-1} there are injective objects in the category ${\rm
Mod}_0S$ and a derived functor must vanish on them.

In this situation semigroups categorical at zero are of special
interest. We recall that a semigroup $S$ with a zero is called {\it
categorical at zero} if for any $x,y,z \in S$ from $xyz=0$ it
follows either $xy=0$ or $yz=0$. For instance, if we join a new
element $0$ to the set of all morphisms of a small category and set
product of morphisms equal to 0 when their composition is not
defined, then the obtained set becomes a semigroups categorical at
zero.

\begin{thrm}{\rm \cite{nov1}}\label{thrm-2-2}
If semigroup $S$ is categorical at zero then $\varepsilon^n$ is an
isomorphism for any $n\ge 0$.
\end{thrm}

This theorem is used in two ways. On the one hand, since, for
instance, completely 0-simple semigroups are categorical at zero
then by Theorem \ref{thrm-2-2} one succeeds to calculate their
0-cohomology, reducing this problem to EM-cohomology \cite{nov3}. On
the other hand, in concrete examples usually it is easier to
calculate 0-cohomology of a given semigroup, and then use it for
finding EM-cohomology of its gown. In more detail we will consider
this question in the following section.

Example \ref{exm-2-3} shows that in general Theorem
\ref{thrm-2-2} is not satisfied.

\section{Calculating EM-cohomology}\label{sec-3}

The free product of semigroups gives a first example of the using
0-cohomo\-logy:

\begin{thrm}\label{thrm-3-1}
Let $S,T$ be semigroups, $S*T$ their free product. Then
$$
H^n(S*T,A)\cong H^n(S,A)\oplus H^n(T,A)
$$
for any $n\ge 2$ and for any $(S*T)$-module $A$.
\end{thrm}

\noindent {\bf Proof.} 0-Direct union $S^0\bigsqcup_0 T^0$ is a
semigroup categorical at zero. Besides, $\widetilde{S^0\bigsqcup_0
T^0}\cong S*T$. So
\begin{eqnarray}
H^n(S*T,A) & \cong & H^n_0(S^0\mbox{$\bigsqcup_0$} T^0,A)\cong
H^n_0(S^0,A)\oplus H^n_0(T^0,A) \nonumber \\
& \cong & H^n(S,A)\oplus H^n(T,A). \nonumber
\end{eqnarray}
Here, the first isomorphism follows from Theorem \ref{thrm-2-2},
second is checked directly, and third follows from Example
\ref{exm-2-2}.

\medskip

{\bf Remark.}  There is an analogue of Theorem \ref{thrm-3-1} for
groups, but its proof is more complicated.

\medskip

In Section \ref{sec-1} the notion of
cohomological dimension was already mentioned. A counter-example to Mitchell conjecture
was obtained by using 0-cohomology. This is the semigroup
$$
S=<a,b,c,d\mid ab=cd>.
$$
Its subset $I=S\setminus\{a,b,c,d,ab\}$ is an ideal and
$\widetilde{S/I}\cong S$. By Corollary \ref{cor-2-2} $H^2(S,A)$ is
embedded into $H^2_0(S/I,A)$.

On the other hand, it is easy to see that $f=\partial\varphi$ for
any 0-cocycle $f\in Z^2_0(S/I,A)$ if we set
$$
\varphi(u)=\left \{
\begin {array}{ll}f(a,b), &\mbox{\rm if}\quad u=a, \\
f(c,d), &\mbox{\rm if}\quad u=c,\\
0 &\mbox{\rm otherwise.}
                  \end {array} \right.
$$
Therefore $H^2(S,A)=H^2_0(S/I,A)=0$.

I give one more example of calculation of EM-cohomology for a pair
of anti-isomorphic semigroups.

Let $p,q\in \mathbb{N}$. Consider the semigroup
$$
T=\langle a,b\mid ab=b, a^p=a^q\rangle.
$$
It is the gown for the semigroup
$$
S=\langle x,y\mid xy=y, x^p=x^q, yx=y^2=0\rangle.
$$
The last is categorical at zero and consists of elements $x^k$
$(k>0)$, $y$ and 0. Its 0-cohomology is easily computed directly, and
thus we get:
\begin{prop}\label{prop-3-1}{\rm \cite{nov4}}
For any $T$-module $A$

1) $H^1(T,A)=H^1_0(S,A)\cong A/\{m\in A\mid am=2m\}$;

2) $H^n(T,A)=H^n_0(S,A)=0$ for $n\ge 0$.
\end{prop}

For the semigroup $T^{\rm op}$ anti-isomorphic to $T$, i.\,e.
$$
T^{\rm op}=\langle a,b\mid ba=b, a^p=a^q\rangle,
$$
the answer is more complicated:

\begin{prop}\label{prop-3-2}{\rm \cite{nov4}}
Let $A$ be an arbitrary $T^{\rm op}$-module, $A_1$ its additive
group, considered as a $T^{\rm op}$-module with trivial
multiplication, $\langle a \rangle$ the subsemigroup generated by
the element $a$. Let homomorphisms $\psi^n$ be induced by embedding
$\langle a \rangle\to T^{\rm op}$. Then the sequence
\begin{eqnarray}
0\to H^0(T^{\rm op},A)\mathop{\to}\limits^{\psi^0} H^0(\langle a
\rangle, A)\to H^0(\langle a \rangle,A_1)\to H^1(T^{\rm
op},A)\mathop{\to}\limits^{\psi^1} ...
\nonumber\\
...\to H^n(T^{\rm op},A)\mathop{\to}\limits^{\psi^n} H^n(\langle a
\rangle,A)\to H^n(\langle a \rangle,A_1)\to ...\nonumber
\end{eqnarray}
is exact.

In particular, if $A_1$ is torsion-free then $H^n(T^{\rm op},A)\cong
H^n(\langle a \rangle,A)$ for $n\ge 1$.
\end{prop}

More general results in calculation of EM-cohomology were
obtained by partial cohomology (see Sec. \ref{sec-8}).

Finally let me formulate an unsolved problem. Theorem \ref{thrm-1-2} shows
that in some cases the cohomology of a semigroup is defined by
the cohomology of its ideal. Generally speaking, this is not the case. When
$I$ is a two-sided ideal, it is desirable to use for
calculation of $H^n(S,A)$ not only the cohomology of the ideal, but also
of the quotient semigroup $S/I$. However, EM-cohomology of the latter
is always trivial, because $S/I$ contains the zero element. So a
question arises: how the group $H^n(S,A)$ depends on $H^n(I,A)$
and $H^n_0(S/I,A)$ (as well as, maybe, on the cohomology groups of
smaller dimension)?

\section{Projective representations}\label{sec-4}

In this section we use the following notation: $S$ is an arbitrary
semigroup (for simplicity we suppose that it contains an identity),
$K$ is a field, $K^{\times}$ its multiplicative group, ${\rm
Mat}_nK$ the multiplicative semigroup of all $(n\times n)$-matrices
over $K$ (we will often delete the subscript $n$). Define an
equivalence $\lambda$ on ${\rm Mat}_nK$: for $A,B\in {\rm Mat}_nK$
put
$$
A\lambda B \Longleftrightarrow  \exists \,c\in K^{\times}:\  A=cB.
$$
Evidently $\lambda$ is a congruence of ${\rm Mat}_nK$. The quotient
semigroup ${\rm PMat}_nK = {\rm Mat}_nK/\lambda$ will be called {\it
a semigroup of projective $(n\times n)$-matrices}.

Let $\Delta:S\to {\rm PMat}_nK$ be a homomorphism and $\alpha: {\rm
Mat}_nK\to {\rm PMat}_nK$ be the canonical homomorphism,
corresponding to congruence $\lambda$. We fix a mapping $\beta: {\rm
PMat}_nK\to {\rm Mat}_nK$, choosing representatives in
$\lambda$-classes. If we denote $\Gamma=\beta\Delta$ then
$\Delta=\alpha\beta\Delta=\alpha\Gamma$. Since $\Delta$ and $\alpha$
are homomorphisms,
$$
\alpha\Gamma(xy)=\Delta(x)\Delta(y)=\alpha(\Gamma(x)\Gamma(y))
$$
\noindent for all $x,y\in S$. Hence $\Gamma(xy)$ and
$\Gamma(x)\Gamma(y)$ vanish simultaneously. So we come to the
following definition:

\begin{defn}\label{defn-4-1}
The mapping $\Gamma: S\to {\rm Mat}_nK$ is called a projective
representation\footnote{One also calls the homomorphism $\Delta$ by a
projective representation.}\label{page-10} of $S$ over $K$ if it
satisfies the following conditions:

1) $\Gamma(xy)=0 \Longleftrightarrow \Gamma(x)\Gamma(y)=0$ for all
$x,y\in S$;

2) there is a partially defined mapping $\rho:S\times S\to
K^{\times}$ such that
\begin{equation}\label{eq-4-1}
{\rm dom\,}\rho= \{(x,y)\mid \Gamma(xy)\ne 0\}
\end{equation}
and
\begin{equation}\label{eq-4-2}
\forall (x,y)\in{\rm dom\,}\rho \quad\Gamma(x)\Gamma(y) =
\Gamma(xy)\rho(x,y).
\end{equation}
The mapping $\rho$ is called a {\it factor set} of  $\Gamma$ and the number $n$
{\it the degree} of $\Gamma$.
\end{defn}

{\bf Remark.} It is easy to see that (\ref{eq-4-2}) remains valid
for all $x,y\in S$ if we extend $\rho$ to a completely defined
mapping setting $\rho(x,y)=0$ for $x, y$ such that $\rho(x,y)$ was
not defined. Hereinafter we will often suppose this.

\medskip

As in the case of projective representations of groups, it is
desirable to have  independent characterization of partially
defined mappings $\rho:S\times S\to K^{\times}$ which can serve as
factor sets for some projective representations of $S$. Applying
(\ref{eq-4-2}) to the equality $\Gamma(x)[\Gamma(y)\Gamma(z)]=
[\Gamma(x)\Gamma(y)]\Gamma(z)$, we get:
\begin{equation}\label{eq-4-3}
\rho(x,y)\rho(xy,z)=\rho(x,yz)\rho(y,z)
\end{equation}
for all $x,y,z\in S$. However, unlike projective representations of
groups, condition (\ref{eq-4-3}) is not sufficient.

\begin{thrm}{\rm \cite{nov2}}\label{thrm-4-1}
A mapping $\rho:S\times S\to K$ is a factor set for a certain
(possible, infinite-dimensional) projective representation of a
monoid $S$ if and only if $\rho$ satisfies (\ref{eq-4-3}) and for all $x,y\in S$
\begin{equation}\label{eq-4-4}
\rho(x,y)=0 \Longleftrightarrow \rho(1,xy)=0.
\end{equation}
\end{thrm}

Similarly to groups, the choice of different representatives of the
$\lambda$-classes leads to an equivalent projective representation.
So we call factor sets $\rho$ and $\sigma$ {\it equivalent}
($\rho\sim \sigma$) if they vanish simultaneously and there exists a
function $\alpha: S\to K^{\times}$ such that for all $x,y\in S$ we have
$$
\rho(x,y)=\alpha(x)\alpha(xy)^{-1}\alpha(y)\sigma(x,y).
$$

Define the product of factor sets $\rho$ and $\sigma$ by pointwise
multiplication: $\rho\sigma(x,y)=\rho(x,y)\sigma(x,y)$. It follows
immediately from Theorem \ref{thrm-4-1} that $\rho\sigma$ is also a
factor set. So the set $m(S)$ of all factor sets is a semigroup and
$\sim$ is its congruence. The quotient semigroup $M(S)=m(S)/\sim$ is
called {\it a Schur multiplier} of $S$.

For groups the Schur multiplier is isomorphic to the group $H^2(G,
K^\times)$ \cite{c-r}. In our situation it is a commutative inverse
semigroup. Consider the construction of the semigroups $M(S)$ and
$m(S)$.

Since $m(S)$ and $M(S)$ are commutative, it follows from Clifford
theorem \cite{c-p} that they are strong semilattices of groups:
$$
m(S)=\bigcup_{\alpha\in b(S)}m_{\alpha}(S),\qquad M(S)=
\bigcup_{\alpha\in B(S)} M_{\alpha}(S),
$$
where $b(S)$ and $B(S)$ are semilattices, $m_{\alpha}(S)$ and
$M_{\alpha}(S)$ are groups. We will call $m_{\alpha}(S)$ and
$M_{\alpha}(S)$ {\it components} of semigroups $m(S)$ and $M(S)$
respectively.

The first step in our consideration is a description of idempotent
factor sets:

\begin{lem}\label{lem-4-1}
There is a bijection $\varepsilon\leftrightarrow I_{\varepsilon}$
between idempotents $\varepsilon\in m(S)$ and ideals of $S$ such
that
$$
\varepsilon(x,y)= \left \{  \begin{array}{ll}
1, & \mbox{\rm if \ } xy\not\in I_\varepsilon, \\ 0, & \mbox{\rm if \ } xy\in I_\varepsilon,
\end{array}  \right.
$$
and
\begin{equation}\label{eq-4-5}
I_{\varepsilon_1\varepsilon_2}=I_{\varepsilon_1}\cup
I_{\varepsilon_2}.
\end{equation}
\end{lem}

We will denote by $Y(S)$ the semilattice of all (two-sided) ideals
of $S$ with respect to union. We consider the empty subset as an
ideal too, i.\,e. $\emptyset \in Y(S)$.

\begin{cor}\label{cor-4-1}
$b(S)\cong B(S)\cong Y(S)$.
\end{cor}

It follows that the ideals $I\in Y(S)$ can serve as
indexes for components of the semigroups $m(S)$ and $M(S)$; thus
$$
m_I(S)m_J(S)\subseteq m_{I\cup J}(S), \qquad M_I(S)M_J(S)\subseteq
M_{I\cup J}(S).
$$

Let $\varepsilon_I$ be the identity of the group $m_I(S)$. Then
$$
\varepsilon_I(x,y)=0\Longleftrightarrow xy\in I.
$$

\begin{lem}\label{lem-4-2}
The group $m_I(S)$ consists of factor sets $\rho$ for which
$$
\rho(x,y)=0\Longleftrightarrow xy\in I.
$$
\end{lem}

Hence groups $m_I(S)$ and $m_0(S/I)$ are isomorphic for
$I\ne\emptyset$. If $I=\emptyset$ we have $m_{\emptyset}(S)\cong
m_0(S^0)$. Certainly, this holds for the multiplier too:

\begin{cor}\label{cor-4-2}
$M_I(S)\cong M_0(S/I)$ if $I\ne\emptyset$, and
$M_{\emptyset}(S)\cong M_0(S^0)$.
\end{cor}

Finally, it is easy to see that $M_0(S)\cong H^2_0(S,K^{\times})$,
and we get the final result:

\begin{thrm}{\rm \cite{nov2}}\label{thrm-4-2}
The Schur multiplier $M(S,K)$ of a semigroup $S$ over a field $K$ is
isomorphic to the semilattice $Y$ of Abelian groups $H^2_0(S/I,
K^{\times})$, where $I\in Y$, and $K^{\times}$ is considered as a
trivial 0-module over $S/I$.
\end{thrm}

Further description of projective representations of semigroups was
carried out in \cite{nov2}; it is similar to description of linear
representations \cite{c-p}.

\section{Brauer monoid}\label{sec-5}

In several articles Haile, Larson and Sweedler
\cite{hai1,hai2,h-l-s,swe}, see also \cite{h-r}, studied so called
strongly primary algebras. Their definition is rather bulky and we
will not need it. Instead of this I cite their description, given in
\cite{hai2}.

Let $K/L$ be a finite Galois extension with the Galois group $G$.
{\it A weak 2-cocycle} \cite{swe} is defined as a mapping
$f:G\times G \to K$ such that for any $\sigma,\tau,\omega \in G$
$$
\sigma[f(\tau,\omega)]f(\sigma\tau,\omega)=f(\sigma,\tau)f(\sigma\tau,\omega),
$$
$$
f(1,\sigma)=f(\sigma,1)=1
$$
(hence weak 2-cocycles can take zero value unlike usual cocycles).

Let $f$ be a weak 2-cocycle. On the set $A$ of formal sums of the
form $\displaystyle{\sum_{\sigma\in G}a_{\sigma}\sigma}$,
$a_{\sigma}\in K$, we define a multiplication by the rule:
$$
a\sigma\cdot b\tau=a\sigma(b)f(\sigma,\tau)\sigma\tau, \qquad
\sigma,\tau\in G,\quad a,b\in K.
$$
Then $A$ becomes an associative algebra. The class of such algebras
coincides with class of strongly primary algebras.

Strongly primary algebras give a generalization of central simple
algebras. In accordance with this Haile, Larson and Sweedler
introduced a notion of a Brauer monoid as generalization of a Brauer
group. For this aim an equivalence of weak 2-cocycles is defined:
$f\sim g$ if there exists a mapping $p:G \to K^{\times}$ such that
$$
g(\sigma,\tau)=f(\sigma,\tau)p(\sigma)p(\tau)(p(\sigma\tau))^{-1}
$$
for any $\sigma,\tau\in G$ under condition $f(\sigma,\tau)\ne 0$.
After factorization by this equivalence the set of weak 2-cocycles
turns into the {\it Brauer monoid} $Br(G,K)$ which is an inverse
semigroup like the Schur multiplier from Sec.\,\ref{sec-4}. More
exactly, denoting by $E$ the semilattice of all idempotents from
$Br(G,K)$ (i.\,e. weak cocycles taking only values 0 and 1), we
obtain:

\begin{thrm} {\rm \cite{h-l-s}}\label{thrm-5-1}
$Br(G,K)$ is a semilattice $E$ of Abelian groups $Br_e(G,K)$, where
$e\in E$ and $Br_e(G,K)$ consists of all weak 2-cocycles which
vanish simultaneously with $e$. In particular, if $e\equiv 1$ then
$Br_e(G,K)\cong H^2(G,K^{\times})$ is the Brauer group.
\end{thrm}

It turns out \cite{k-n,nov10} that this construction is reduced to
0-cohomology. Let $e\in E$. Join an extra zero 0 to $G$ and define a
new operation on $G^0$ :
$$
x\circ y=\left \{ \begin {array}{ll}xy, & \mbox{\rm if}\ e(x,y)=1, \\
0, & \mbox{\rm if}\ e(x,y)= 0
\end {array} \right.
$$
and besides, $x\circ 0=0\circ x=0$. With this operation $G^0$ is a
semigroup which we will denote by $G_e$. Conversely, call by a {\it
modification} $G(\circ)$ of the group $G$ a monoid on $G^0$ with an
operation $\circ$ such that $x\circ y$ is either $xy$ or 0, and
moreover $0\circ x=x\circ 0=0$. It is easy to see that there is a
bijective correspondence between idempotent weak 2-cocycles and
modifications of $G$. The group $K^{\times}$ turns into a 0-module
over $G_e$, $Br_e(G,K)\cong H^2_0(G_e,K^{\times})$ and Theorem
\ref{thrm-5-1} changes into the following statement:

\begin{thrm}\label{thrm-5-2}
$Br(G,K)$ is a semilattice of Abelian groups
$H^2_0(G(\circ),K^{\times})$, where $G(\circ)$ runs the set of all
modifications of the group $G$.
\end{thrm}

It is shown in \cite{nov10} that in this problem 0-cohomology is
used essentially: to describe some properties of Brauer monoid one
has to use 0-cohomology of other (different from modifications)
semigroups.

Thus study of the Brauer monoid is reduced to description of
modifications of the group and their 0-cohomology. However, it is
necessary to note that the study of modifications is a difficult
combinatorial problem. In general case for a finite group $G$ (only
such groups are considered in Haile--Larson--Sweedler theory) we can
only confirm that each modification is an union of the subgroup of
its invertible elements and a nilpotent ideal. Besides,
modifications are 0-cancellative (if $ax=bx\ne 0$ or $xa=xb\ne 0$
then $a=b$). Some examples of modification were considered in
\cite{nov11} and \cite{nov14}.

\section{Partial representations of groups}\label{sec-6}

Results of this section were received when I worked in Sa\~o Paulo,
Brasil, thanks to the foundation FAPESR. They were announced on XVIII
Brazilian Algebra Meeting \cite{d-n} and are preparing for
publication.

In connection with studying $C^*$-algebras so called partial linear
representations of groups appeared \cite{exe,q-r}. It is naturally
to ask: how do partial {\bf projective} representations of groups
look like? It turned out that here 0-cohomology appears too. We start
with necessary definitions from \cite{exe}.

\begin{defn}\label{defn-6-1}
A mapping $\varphi:G\to S$ from a group $G$ into a semigroup $S$ is
called a partial homomorphism if for all $x,y\in G$
\begin{eqnarray}
\varphi(x^{-1})\varphi(x)\varphi(y)&=&\varphi(x^{-1})\varphi(xy)\nonumber\\
\varphi(x)\varphi(y)\varphi(y^{-1})&=&\varphi(xy)\varphi(y^{-1})\nonumber\\
\varphi(x)\varphi(e)&=&\varphi(x)\nonumber
\end{eqnarray}
(these equalities imply $\varphi(e)\varphi(x)=\varphi(x)$).

In particular, a partial linear representation (PLR) over a field
$K$ is a partial homomorphism into the matrix semigroup, $\Delta:
G\to{\rm Mat}_nK$.
\end{defn}

R.\,Exel introduced a monoid $\Sigma(G)$ which plays a special
role here. It is generated by symbols $[x]$ ($x\in G$) with defining
relations
\begin{eqnarray}
[x^{-1}][x][y]&=&[x^{-1}][xy]\nonumber\\
{}[x][y][y^{-1}]&=&[xy][y^{-1}]\nonumber\\
{}[x][e]&=&[x]\nonumber
\end{eqnarray}
(these equalities imply $[e][x]=[x]$).

$\Sigma(G)$ possesses the following universal property:
\begin{enumerate}
\item The mapping $f:G\to \Sigma(G)$, $f(x)=[x]$, is a partial homomorphism.

\item For any semigroup $S$ and any partial homomorphism $\varphi:G\to S$
there exists an unique (usual) homomorphism $\tilde{\varphi}:
\Sigma(G)\to S$ such that $\varphi=\tilde{\varphi}f$.
\end{enumerate}

Due to this property  the study of PLR's of groups is equivalent to
the study of linear representations of its Exel monoid.

It is naturally to define (and to study) partial projective
representations of groups by means of usual projective
representations of $\Sigma(G)$: we will call the partial
homomorphism $\Delta : G \to {\rm PMat}K$ by a {\it partial
projective representation (PPR)} of $G$ (cf. the footnote on
p.\,\pageref{page-10}). We get the diagram

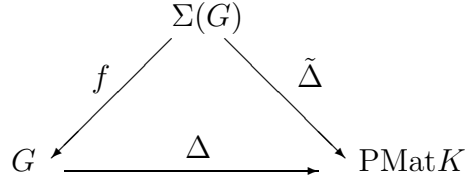
\begin{picture}(100,80)
\put(200,65){$\Sigma(G)$} \put(140,10){$G$}
 \put(270,10){${\rm PMat}K$}
\put(160,10){\vector(1,0){95}} \put(205,15){$\Delta$}
\put(200,60){\vector(-1,-1){45}} \put(170,42){$f$}
\put(220,60){\vector(1,-1){45}} \put(247,42){$\tilde{\Delta}$}
\end{picture}

\noindent where $\tilde{\Delta}$ is a projective representation of
$\Sigma(G)$.

However we will see below that PPR's are not reduced to projective
representations of semigroups unlike linear ones.

The first step in study of PPR's is a translation of their
definition into the language of usual matrices:

\begin{thrm}\label{thêü-6-1}
A mapping $\Gamma: G\to {\rm Mat}K$ is PPR of $G$ if and only if the
following conditions hold:

\medskip

1) for all $x,y\in G$
$$
\Gamma (x^{-1})\Gamma (xy) = 0 \Longleftrightarrow
\Gamma(x)\Gamma(y) = 0 \Longleftrightarrow \Gamma(xy)\Gamma(y^{-1})
= 0;
$$

2) there is a mapping $\sigma:G\times G\to K$ such that
$$
\Gamma(x)\Gamma(y)=0 \Longleftrightarrow \sigma(x,y)=0
$$
and
\begin{eqnarray}
\Gamma (x^{-1})\Gamma (x)\Gamma (y) &=&  \Gamma (x^{-1})\Gamma
(xy) \sigma (x,y),\nonumber\\
\Gamma (x)\Gamma (y)\Gamma (y^{-1}) &=&  \Gamma (xy)\Gamma (y^{-1})
\sigma(x,y).\nonumber
\end{eqnarray}
\end{thrm}

Note that this theorem gives another definition of PPR independent of $\Sigma(G)$.

We call $\sigma$ {\it a factor set} of $\Gamma$ and define a product
of factor sets as above. However, it is not evident that this
product is also a factor set. To prove this, one have to use the Exel monoid again:

\begin{prop}\label{prop-6-1}
Let $\sigma:G\times G\to K$ be a mapping for which there is a factor
set $\rho$ of the semigroup $\Sigma(G)$, such that:
\begin{eqnarray}
&1)&\ \forall\, x,y\in G\quad \sigma(x,y)=0 \Longleftrightarrow
\rho([x],[y])=0;\nonumber\\
&2)&\ \sigma(x,y)\ne 0 \Longrightarrow \sigma(x,y)=
\frac{\rho([x],[y])\,\rho([x^{-1}],[x][y])}{\rho([x^{-1}],[xy])}.\nonumber
\end{eqnarray}
Then $\sigma$ is a factor set of some PPR of $G$.

The converse also holds.
\end{prop}

Now the desired result about product comes directly. Moreover:
\begin{cor}\label{cor-6-1}
The factor sets of $G$ form a commutative inverse semigroup $Pm(G)$.
\end{cor}

Again, as in Sec.\,\ref{sec-4}, we define an equivalence of factor
sets and call the respective quotient semigroup by {\it the Schur
multiplier} $PM(G)$. The Schur multiplier is also a commutative
inverse semigroup. However, the Schur multipliers of $G$ and of
$\Sigma(G)$ are different: one can only confirm that $PM(G)$ is an
image of $M(\Sigma(G))$ (Prop.\,\ref{prop-6-1}).

At present the main problem in description of a Schur multiplier is
to find such a definition of factor sets which would be independent on
both PPR and $\Sigma(G)$. We recall that for usual projective
representations of groups such a definition is the cohomological
equation (\ref{eq-4-3}). In the case of semigroups (more exactly,
monoids) it is necessary to add the condition (\ref{eq-4-4}).
Unfortunately, a cohomological equation is not true for PPR's. At
most we can assert
\begin{prop}\label{prop-6-2}
Let $\Gamma$ be a PPR with a factor set $\sigma$. Then
$$
\forall x,y,z\in G\quad\Gamma(x)\Gamma(y)\Gamma(z)\ne 0
\Longrightarrow \sigma(x,y)\sigma(xy,z)=\sigma(x,yz)\sigma(y,z).
$$
\end{prop}

Later I will give an example where the cohomological equation does
not hold.

Now we consider the structure of a Schur multiplier. Certainly, it
is a semilattice of their subgroups. So first of all it is necessary
to describe its idempotents:
\begin{thrm}\label{thrm-6-2}
Let $\sigma:G\times G\to K$ be a mapping taking only values $0$ and
$1$ and $\sigma(1,1)=1$. Then $\sigma$ is a factor set if and only
if
\begin{equation}\label{eq-6-1}
\forall\,x,y\in G \quad \sigma(x,y)=1 \Longrightarrow
\sigma(xy,y^{-1})=\sigma(y^{-1},x^{-1})=\sigma(x,1)=1.
\end{equation}
\end{thrm}

\begin{exm}\label{exm-6-1}
Let $G=\langle a,b,c\rangle$ be an elementary Abelian group of order
8 with generators $a,b,c$. Denote $H=\langle b,c\rangle$,
$F=(H\setminus 1)\times (H\setminus 1)\setminus \nabla$ where
$\nabla$ is a diagonal of the Cartesian square $H\times H$. Set
$$
\sigma(x,y)= \left\{\begin{array}{cc}
1 &\mbox{\rm if\ } (x,y)\not\in F,\\ 0 &\mbox{\rm if\ } (x,y)\in F.
\end{array}
\right.
$$
It is easy to check by Theorem \ref{thrm-6-2} that $\sigma$ is a
factor set. However,
$$
\sigma(b,a)\sigma(ba,ac)=1\ne 0=\sigma(b,c)\sigma(a,ac)
$$
and the cohomological equation does not hold for $x=b$, $y=a$,
$z=ac$.
\end{exm}

One can place Theorem \ref{thrm-6-2} into a more general form.
Consider an abstract semigroup $\cal T$ generated by elements
$\alpha,\beta,\gamma$ with defining relations
$$
\left.
\begin{array}{cl}
&\alpha^2 =\beta^2=1,\quad (\alpha \beta)^2=1\\
&\gamma^2 =1,\quad \alpha\gamma=\gamma,\quad
\gamma\alpha\beta\gamma=\gamma\beta\alpha\beta,\quad
\gamma\beta\gamma=0
\end{array}
\right\}
$$
For any group $G$ the semigroup $\cal T$ acts on $G\times G$ as
follows:
\begin{eqnarray}
\alpha : (x,y)&\longrightarrow& (xy,y^{-1})\nonumber\\
\beta : (x,y)&\longrightarrow& (y^{-1},x^{-1})\nonumber\\
\gamma : (x,y)&\longrightarrow& (x,1)\nonumber
\end{eqnarray}

Thus, $G\times G$ turns into $\cal T$-set and Theorem \ref{cor-6-1}
takes the form:

\begin{cor}\label{cor-6-2}
An idempotent mapping $\sigma:G\times G\to K$ such that
$\sigma(1,1)=1$ is a factor set if and only if ${\rm
supp\,}\sigma=\{(x,y)\mid \sigma(x,y)\ne 0\}$ is a $\cal T$-subset
in $G\times G$.
\end{cor}

Now from Corollaries \ref{cor-6-1} and \ref{cor-6-2} we get:

\begin{thrm}\label{thrm-6-3}
The Schur multiplier is a semilattice of Abelian groups
$$
Pm(G)=\bigcup_{X\in C(G)}Pm_X(G),\qquad PM(G)= \bigcup_{X\in C(G)}
PM_X(G),
$$
where $C(G)$ is a semilattice of $\cal T$-subsets in $G\times G$
with respect to intersection.
\end{thrm}

Let me say a few words about the semigroup $\cal T$. It plays a remarkable
role: for any group $G$ it gives a description of idempotent factor
sets. Since each PLR is a PPR with an idempotent factor set, we
obtain, in particular, some classification of all PLR's of a group
$G$. So $\cal T$  merits to be considered more thoroughly. Here are more details on its structure.

First of all, its order is 25. The elements $\alpha$ and $\beta$
generate in $\cal T$ a subgroup $H=\langle \alpha,\beta\rangle$,
isomorphic to the symmetric group $S_3$. The complement $U= {\cal T}
\setminus H$ is an ideal.

One can prove that $U$ is a completely 0-simple semigroup. In the
standard notation of Theory of Semigroups \cite{c-p} it can be
written as $U=M^0(D;I,\Lambda;P)$, where $I=\Lambda=\{1,2,3\}$, $D$
is a group of the order 2 and $P$ is a $(3\times
3)$-sandwich-matrix,
$$
P=\left( \begin{array}{ccc}
1&1&0\\
1&0&1\\
0&1&1
\end{array}
\right).
$$

\section{Cohomology of small categories}\label{sec-7}

In further study of properties of 0-cohomology some difficulties
arise because, as I mentioned already, 0-cohomology is not a derived
functor in the Abelian category where it is built (see
Ex.\,\ref{exm-2-3}).

So a question appears: is it possible to extend the category of
0-modules so that 0-cohomology becomes a derived functor? One of the
useful ways is to pass to bimodules. However, in our situation it
does not help as one can see from the following example.

We call an Abelian group $A$ by a {\it $0$-bimodule} over $S$ if $A$
is right and left 0-module, and besides, $(sa)t=s(at)$ for any
$s,t\in S\setminus 0$, $a\in A$. 0-Cohomology of $S$ with values in
the category of 0-bimodules are defined similarly to 0-cohomology on
0-modules. Denote it by $HH^n_0(S,A)$.

\begin{exm}\label{exm-7-1}
Let $S=\{u,v,w,0\}$ be a commutative semigroup with multiplication
$u^2=v^2=uv=w, \ uw=vw=0$, $M$ a $0$-bimodule over $S$. Then
$HH^2_0(S,M)\ne 0$ for $M\ne 0$.
\end{exm}

As in Sec.\,\ref{sec-2}, this example shows that in the category of
0-bimodules the cohomology functor $HH^n_0$ is not a derived
functor. This is the reason why we use the category  $\NatS$ (which
is defined below). Our construction is a generalization of the
theory of cohomology for small categories from \cite{b-w}.

\medskip

As in Sec.\,\ref{sec-4}, we
suppose for simplicity that $S$ is a monoid with a zero. Call by the
{\it category of factorizations in} $S$ the category ${\cal F}acS$
whose objects are all elements from $S\setminus0$, and the set of
morphisms ${\rm Mor}(a,b)$ consists of all triples
$(\alpha,a,\beta)$  ($\alpha,\beta \in S$) such that $\alpha a\beta
= b$ (we will denote $(\alpha,a,\beta)$ by $(\alpha,\beta)$ if this
does not lead to confusion). The composition of morphisms is defined
by the rule $(\alpha',\beta')(\alpha,\beta)=
(\alpha'\alpha,\beta\beta')$;  hence we have $(\alpha,\beta)=
(\alpha,1)(1,\beta)=(1,\beta)(\alpha,1)$.

A {\it natural system on} $S$ is a functor ${\bf D}:{\cal F}acS \to
{\cal A}b$. The category $\NatS = {\cal A}b ^{{\cal F}acS}$ of such
functors is an Abelian category with enough projectives and
injectives \cite{gro}. Denote the value of ${\bf D}$ on an object
$a\in {\rm Ob} {\cal F}acS$ by ${\bf D}_a$.  If we denote
$\alpha_*={\bf D}(\alpha,1)$ and $\beta^*={\bf D}(1,\beta)$ then
${\bf D}(\alpha,\beta)=\alpha_*\beta^*$ for any morphism
$(\alpha,\beta)$.

\begin{exm}\label{exm-7-2}
Each 0-module $A$ can be considered as a functor ${\bf A}$ from
$\NatS$, defined as follows: ${\bf A}_s =A$ for any $s\in S\setminus
0$ and $\alpha_*\beta^*a =\alpha a$ for all $\alpha,\beta\in S$,
$a\in A$.
\end{exm}

\begin{exm}\label{exm-7-3}
Consider a functor ${\bf Z}$ which assigns to each object $a\in
S\setminus 0$ the infinite cyclic group ${\bf Z}_a$ generated by a
symbol $[a]$; to each morphism $(\alpha,\beta):s\to t$ it assigns a
homomorphism of the groups ${\bf Z}(\alpha,\beta):{\bf Z}_a\to {\bf
Z}_b$ which takes $[a]$ to $[b]$. It is a natural system, which is
called trivial.
\end{exm}

For a given natural number $n$ denote by $Ner_n S$ the set of all
$n$-tuples $(a_1,\ldots,a_n)$, $a_i\in S$, such that $a_1\cdots
a_n\ne 0$ (a {\it nerve} of $S$). For $n=0$ we set $Ner_0 S=\{1\}$.
A mapping, defined on the nerve and assigning to each
$a=(a_1,\ldots,a_n)$ an element from ${\bf D}_{a_1\cdots a_n}$, is
called an {\it $n$-dimensional cochain}. The set of all
$n$-dimensional cochains is an Abelian group $C^n(S,{\bf D})$ with
respect to the pointwise addition. Set $C^0(S,{\bf D})={\bf D}_1$.

Define a {\it coboundary homomorphism} $\Delta^n:C^n(S,{\bf D})\to
C^{n+1}(S,{\bf D})$ for $n\geq 1$ by the formula
\begin{eqnarray}
&&(\Delta^n f)(a_1,\ldots,a_{n+1})=a_1{}_* f(a_2,\ldots,a_{n+1})
\nonumber \\
&&+\sum_{i=1}^{n}{(-1)^i f(a_1,\ldots,a_i a_{i+1},\ldots
,a_{n+1})}+(-1)^{n+1}a_{n+1}^*f(a_1,\ldots,a_n).\nonumber
\end{eqnarray}

For $n=0$ we set $(\Delta^0 f)(x)=x_*f-x^*f$ where $f\in {\bf D}_1$,
$x\in S\setminus 0$. One can check directly that
$\Delta^n\Delta^{n-1}=0$. Cohomology groups of the complex
$\{C^n(S,{\bf D}),\Delta^n \}_{n\geq0}$ are denoted by $H^n(S,{\bf
D})$.

0-Cohomology of a monoid is a special case of this construction.
Namely, $H_0^n(S,A)\cong H^n(S,{\bf A})$, where ${\bf A}$ is a
functor defined in Ex.\,\ref{exm-7-2}.

Since $\NatS$ has enough projectives and injectives there exist
derived functors ${\rm Ext}^n_\NatS ({\bf Z},\rule{6pt}{.8pt})$.

\begin{thrm}\label{thrm-7-1}{\rm \cite{k-n1,k-n2}}
For any monoid $S$ with zero
$$
H^n(S,\rule{6pt}{.8pt})\cong{\rm Ext}_{\NatS}^n({\bf
Z},\rule{6pt}{.8pt}).
$$
\end{thrm}

To prove this statement a projective resolution for ${\bf Z}$ is
built in the following way.

For every $n\geq 0$ we define a natural system ${\bf B}_n:{\cal F}ac
S\to {\cal A}b$. Its value on an object $a\in S\setminus 0$ is a
free Abelian group ${\bf B}_n(a)$ generated by the set of symbols
$[a_0,\ldots,a_{n+1}]$ such that $a_0\cdots a_{n+1}=a$. To each
morphism $(\alpha,\beta)$ we assign a homomorphism of groups by the
formula
$$
{\bf B}_n(\alpha, \beta):\ [a_0,\ldots,a_{n+1}] \to [\alpha
a_0,\ldots,a_{n+1}\beta].
$$
These functors constitute a chain complex $\{{\bf B}_n,
\partial_n\}_{n\geq 0}$, where natural transformations
$\partial_n:{\bf B}_n\toop {\bf B}_{n-1} \ \ (n\geq 1)$ are given by
the homomorphisms
$$
(\partial_n)_a:{\bf B}_n(a)\to{\bf B}_{n-1}(a),
$$
$$
(\partial_n)_a[a_0,\ldots,a_{n+1}] =
 \sum_{i=0}^{n}{(-1)^i[a_0,\ldots,a_ia_{i+1},\ldots,a_{n+1}].}
$$

The natural systems ${\bf B}_n$ are projective objects in $\NatS$
and the complex $\{{\bf B}_n,\partial_n\}_{n\geq 0}$ is a projective
resolution of the natural system ${\bf Z}$.

Now one can establish an isomorphism between the complexes
$$
\{C^n(S,{\bf D}),\Delta^n\}_{n\geq 0} \ \ {\rm and} \ \ \{{\rm
Hom}_{\NatS}({\bf B}_n,{\bf D}),\partial^n\}_{n\geq 0}.
$$

\medskip

Our construction differs from Baues cohomology theory for monoids
\cite{b-w} in the initial stage only. Namely, in \cite{b-w} a monoid
$S$ is regarded as a category with a single object. At the same time
the Baues category of factorizations in $S$ is equal to ${\cal F}ac
S^0$ where $S^0$ is a semigroup with a joined zero. Therefore the
Baues cohomology groups of $S$ and the cohomology groups of $S^0$ in
our sense are the same. However if $S$ possesses a zero element then
the category ${\cal F}acS$ and Baues one are not equivalent and we
obtain the different cohomology groups. The construction of this
section is a generalization simultaneously both Baues' and
0-cohomology.

\medskip

In conclusion I give an example of using obtained results.

It is well-known that in many algebraic theories cohomological
dimension of free objects is 1. In the category of monoids with zero
a free object is a free monoid with a joined zero. However in this
category the class of objects having cohomological dimension 1 is
essentially greater.

Call every quotient monoid of a free monoid by its ideal {\it
$0$-free}. Free monoids with a joined zero are also considered as
0-free. Let $S$ be a semigroup with a zero. We call the least $n$
such that $H^{n+1}_0(S,A)=0$ for any $0$-module $A$, by {\it
$0$-cohomological dimension} (0-{\rm cd}$S$) of $S$.

\begin{thrm}\label{thrm-7-2}
$0${\rm -cd}$M\leq 1$ for any $0$-free monoid $M$.
\end{thrm}

From this theorem it follows an interesting

\begin{cor}\label{cor-7-1}
Any projective representation of a $0$-free monoid is linearized
(i.\,e. is equivalent to a linear one).
\end{cor}

In connection with Theorem \ref{thrm-7-2} a question arises, an
answer to which is unknown to me: is a 0-cancellative monoid
(Sec.\,\ref{sec-6}), having 0-cohomological dimension 1, 0-free?

\section{Concluding remarks}\label{sec-8}

General properties of 0-cohomology are studied weakly. Here the same
difficulties (and even greater) appear as for EM-cohomology of
semigroups. Certain hopes are given by Theorem \ref{thrm-7-1},
showing that 0-cohomology can be continued to a derived functor.
However this continuation turns out to be too vast. So it remains
actual to find a category, smaller than ${\cal F}ac$,
in which 0-cohomology would be a derived functor.

0-Cohomology appear in other problems too. In an article by Clark
\cite{cla} it was applied to semigroups of matrix units and
algebras generated by them. Similar situation often occurs in Ring Theory.
A quotient algebra of a semigroup algebra by its ideal,
generated by the zero of the semigroup, is called a {\it contracted}
one (in other words, the zero of the semigroup is identified with
the zero of the algebra). For instance, the well-known theorem of
Bautista--Gabriel--Roiter--Salmeron \cite{b-g-r-s} confirms that
every algebra of the final type is contracted semigroup one.

A natural question arises: how the Hochschild cohomology of
contracted algebras is connected with the cohomology of semigroups generating them?
Since it is supposed that the semigroup contains a zero, then, of
course, for study of this question it is necessary to use
0-cohomology. Such approach could be useful for incidence algebras
of simplicial complexes as well (cf. \cite{g-s}).

In conclusion let me mention a generalization of 0-cohomology. In a
semigroup $S$ (not necessary containing 0) let us fix certain
generating subset $W\subset S$ instead of $S\setminus 0$ and call
mappings $W\to A$ by {\it $1$-dimensional $W$-coboundaries}. Using this
one can construct certain partial $n$-place mappings
$S\times\ldots\times S\to A$ (and call them {\it $n$-dimensional
W-cochains}) so that coboundary homomorphisms $\partial^n$ would be
well defined. I have called the obtained objects by {\it partial
cohomology} (they have no concern with partial representations from
Sec.\,\ref{sec-6}) and considered it in \cite{nov7} and \cite{nov8}.
Partial cohomology turned out to be useful for calculation of
EM-cohomology (as in Sec.\,\ref{sec-3}), however I have not meet
other applications of them.


\begin{thebibliography}{111}

\bibitem{a-r}
Adams W.\,W., Rieffel M.\,A.
        {\it Adjoint functors and derived functors with
        an application to the cohomology of semigroups.}
        J. Algebra, {\bf 7}(1967), N1, 25-34.

\bibitem{b-w}
Baues H.-J., Wirshing G.
        {\it Cohomology of small categories.}
        J.\,Pure Appl.\,Algebra, {\bf 38}(1985), N2/3, 187-211.

\bibitem{b-g-r-s}
Bautista R., Gabriel P., Roiter A.\,V., Salmeron L. {\it
Representation-finite algebras and multiplicative bases.} Invent.
Math., {\bf 81}(1985), N2, 217-286.

\bibitem{bro}
Brown K.\,S. {\it Cohomology of Groups}. Springer,1982.

\bibitem{c-e}
Cartan H., Eilenberg S. {\it Homological Algebra.} Princeton, 1956.

\bibitem{cla}
Clark W.\,E.
        {\it Cohomology of semigroups via topology with an
        application to semigroup algebras.}
        Commun. Algebra, {\bf 4}(1976), 979-997.

\bibitem{c-p}
Clifford A.\,H., Preston G.\,B. {\it Algebraic Theory of
Semigroups.} Amer. Math. Soc., Providence, 1964.

\bibitem{c-r}
Curtis C.\,W., Reiner I. {\it Representation Theory of Finite Groups
and Associative Algebras.} Willey \& Sons, N.-Y., L., 1967.

\bibitem{d-n}
Dokuchaev M., Novikov B. {\it Projective representations and Exel's
theory.} In: XVIII Escola de Algebra (abstracts), Campinas, July
19-23, 2004, 31-32.

\bibitem{exe}
Exel R. {\it Partial actions of groups and actions of semigroups.}
Proc. Amer. Math. Soc., {\bf 126}(1998), N12, 3481-3494.

\bibitem{g-s}
Gerstenhaber M., Schack S.\,D. {\it Simplicial cohomology is
Hochschild cohomology.} J. Pure Appl. Algebra {\bf 30}(1983), N2,
143-156.

\bibitem{gri}
Grillet P.\,A.
        {\it Commutative semigroup cohomology.}
        Semigroup Forum, {\bf 43}(1991), N2, 247-252.

\bibitem{gro}
Grothendieck A. {\it Sur quelques points d'alg\`ebre homologique.}
Tohoku Math. J., {\bf 9}(1957), N2, 119-183; N3, 185-221.

\bibitem{hai1}
Haile D.\,E.
        {\it On crossed product algebras arising from weak cocycles.}
        J. Algebra, {\bf 74}(1982), 270-279.

\bibitem{hai2}
Haile D.\,E.
        {\it The Brauer monoid of a field.}
        J. Algebra, {\bf 81}(1983), N2, 521-539.

\bibitem{h-l-s}
Haile D.\,E., Larson R.\,G., Sweedler M.\,E.
        {\it A new invariant for ${\bf C}$ over ${\bf R}$: almost
        invertible cohomology theory and the classification of idempotent
        cohomology classes and algebras by partially ordered sets
        with Galois group action.}
        Amer. J. Math., {\bf 105}(1983), N3, 689-814.

\bibitem{h-r}
Haile D., Rowen L. {\it Weakly Azumaya algebras.} J. Algebra  {\bf
250}(2002), 134-177.

\bibitem{han}
Hancock V.\,R.
        {\it Commutative Schreier semigroup extensions of a group.}
        Acta Sci. Math., {\bf 25}(1964), N2, 129-134.

\bibitem{h-s}
Hilton P.\,J., Stammbach U. {\it A Course in Homological Algebra.}
Springer, 1971.

\bibitem{k-n}
Kirichenko V.\,V., Novikov B.\,V. {\it On the Brauer monoid for
finite fields.} Finite fields and applications (Augsburg, 1999).
Springer, Berlin, 2001, 313-318.

\bibitem{k-n1}
Kostin A.\,A., Novikov B.\,V.
        {\it Semigroup cohomology as a derived functor.}
        Filomat, {\bf 15}(2001), 17-23.

\bibitem{k-n2}
Kostin A.\,A., Novikov B.\,V. {\it Cohomology of semigroups and
small categories.} In: Algebraic Structures and Their Applications.
Proc. of Ukrain. Math. Congress, Kyiv, 2002, 69-75 (Russian).

\bibitem{lau}
Lausch H.
        {\it Cohomology of inverse semigroups.}
        J. Algebra, {\bf 35}(1975), N1-3, 273-303.

\bibitem{lee}
Leech J.\,E. {\it Two papers: $H$-coextensions of monoids and The
structure of a band of groups.} Mem. Amer. Math. Soc., {\bf
157}(1975).

\bibitem{log}
Loganathan M.
        {\it Cohomology and extensions of regular semigroups.}
        J.~ Austral. Math. Soc., ser. A, {\bf 35}(1983), N2,
        178-193.

\bibitem{mit}
Mitchell B.
    {\it On the dimension of objects and categories. {\rm I.} Monoids.}
    J. Algebra., {\bf 9}(1968), N3, 314-340.

\bibitem{nic}
Nico W.\,R.
        {\it On the cohomology of finite semigroups.}
        J. Algebra, {\bf 11}(1969), N4, 598-612.

\bibitem{nov1}
Novikov B.\,V. {\it $0$-cohomology of semigroups.} In: Theoretical
and applied questions of differential equations and algebra.
"Naukova Dumka", Kiev, 1978, 185-188 (Russian).

\bibitem{nov2}
Novikov B.\,V. {\it Projective representations of semigroups.} Dokl.
Akad. Nauk Ukrain. SSR, Ser. A, N6, 1979, 474-478 (Russian).

\bibitem{nov3}
Novikov B.\,V. {\it $0$-cohomology of completely $0$-simple
semigroups.} Vestnik Kharkov. Gos. Univ. N221, 1981, 80-85
(Russian).

\bibitem{nov4}
Novikov B.\,V. {\it On computing of $0$-cohomology of some
semigroups.} Vestnik Kharkov. Gos. Univ. N221, 1981, 96 (Russian).

\bibitem{nov5}
Novikov B.\,V. {\it A counterexample to a conjecture of Mitchell.}
Trudy Tbiliss. Mat. Inst. Razmadze Akad. Nauk Gruzin. SSR, {\bf 70}
(1982), 52-55 (Russian).

\bibitem{nov6}
Novikov B.\,V. {\it Defining relations and $0$-modules over a
semigroup.} In: Theory of semigroups and its applications. Saratov.
Gos. Univ., Saratov, 1983, 94-99 (Russian).

\bibitem{nov7}
Novikov B.\,V.
        {\it On partial cohomologies of semigroups.}
        Semigroup Forum, {\bf 28}(1984), N1-3, 355-364.

\bibitem{nov8}
Novikov B.\,V. {\it Partial cohomology of semigroups and its
applications.} Izv. Vyssh. Uchebn. Zaved. Mat., 1988, N11, 25-32
(Russian); translation in Soviet Math. (Iz. VUZ) {\bf 32}(1988),
N11, 38-48.

\bibitem{nov9}
Novikov B.\,V. {\it Commutative semigroups with cancellation of
dimension $1$.} Mat. Zametki, {\bf 48}(1990), N1, 148-149 (Russian).

\bibitem{nov10}
Novikov B.\,V.{\it The Brauer monoid.} Matem. zametki, {\bf
57}(1995), N4, 633-636 (Russian). Translated in: Math. Notes {\bf
57}(1995), N3-4, 440-442.

\bibitem{nov11}
Novikov B.\,V.
        {\it On modifications of the Galois group.}
        Filomat, {\bf 9}(1995), N3, 867-872.

\bibitem{nov12}
Novikov B.\,V. {\it Semigroups of cohomological dimension 1.}
J.\,Algebra, {\bf 204}(1998), 386-393.

\bibitem{nov13}
Novikov B.\,V.{\it The cohomology of semigroups: a survey.} Fundam.
Prikl. Mat. {\bf 7}(2001), N1, 1-18 (Russian).

\bibitem{nov14}
Novikov B.\,V. {\it Semigroup cohomology and applications.} In:
Algebra -- Representation Theory (ed. K.\,W.\,Roggenkamp and M.\,\c
Stef\u anescu), Kluwer, 2001, 219-234.

\bibitem{pac}
Pachkoriya, A.\,M. {\it Cohomology of monoids with coefficients in
semimodules.} Soobshch. Akad. Nauk Gruzin. SSR, {\bf 86}(1977), N3,
546-548 (Russian).

\bibitem{q-r}
Quigg J.\,C., Raeburn I. {\it Characterizations of crossed products
by partial actions.} J. Operator Theory, {\bf 37}(1997), 311-340.

\bibitem{str}
Strecker R.
        {\it Uber kommutative Schreiersche Halbgruppenerweiterungen.}
        Acta Math. Acad. Sci. Hung., {\bf 28}(1972), N1-2, 33-44.

\bibitem{swe}
Sweedler M.\,E.
        {\it Weak cohomology.} Contemp. Math., {\bf 13}(1982), 109-119.

\bibitem{wel}
Wells Ch.
        {\it Extension theories for monoids.} Semigroup Forum, {\bf 16}(1978), N1, 13-35.

\end{thebibliography}
\end{document}